
\magnification=\magstephalf 


 \font\smcaps=cmcsc10  
 \font\smrm=cmr7   


  \def\beginsmallparskip{\bgroup \parskip=\smallskipamount}
  \def\endsmallparskip{\par\egroup}
  \def\noparindent{\parindent=0pt}
  
  \def\doubleparindent{\parindent=40pt}


   \def\heading#1{\bgroup\noparindent\smallskip \smcaps #1 \egroup}


   \def\theorem#1{\bgroup\noparindent\smallskip\smcaps #1 \par\egroup}


   \def\card{\hbox{\rm card}}  


   \def\refno#1{\bgroup \bf #1\egroup}
   \def\cite#1#2{\bgroup \bf [#1\egroup, #2\bgroup \bf ]\egroup}


\ 
\bigskip
\centerline{\bf New Results on Primes from an Old Proof of Euler's}
\bigskip
\centerline{\rm by}
\medskip
\centerline{Charles W. Neville}
\smallskip
\centerline{CWN Research}
\centerline{55 Maplewood Ave.}
\centerline{West Hartford, CT \ \ 06119, U{.}S{.}A{.}}
\centerline{cwneville@cwnresearch.com}
\bigskip
\centerline{June, 2002}
\centerline{Revision 1, April, 2003}

\vfill\eject


\beginsmallparskip

\heading{Introduction.}

{\noparindent
In 1737 Leonard Euler gave what we often {\sl now\/} think of as a new proof, based on infinite series, of Euclid's theorem that there are infinitely many prime numbers.  This work was noteworthy in many ways.  It marked the first use of the methods of analysis in number theory, and it marked the first appearance of what later came to be known as the Riemann $\zeta$ function in mathematics.\par}

In actual reality, Euler never presented his work as a proof of Euclid's theorem, though that conclusion is clearly implicit in what he did.  Rather, Euler derived such striking new results (for the time) as $\sum 1/p$ diverges, where $p$ ranges over the prime numbers.

   Briefly, Euler considered the possibly infinite product, $\prod  1/\left( 1 - p^{-1} \right)$, where the index $p$ runs over all primes.  He expanded the product to obtain the divergent infinite series $\sum_{n=1}^\infty  1/n$, concluded the infinite product was also divergent, and from this concluded that the infinite series $\sum  1/p$ also diverges.  Well, if there were only finitely many primes $p$, the series would have to converge because finite series are ALWAYS CONVERGENT.  So we, and Euler had he chosen to do so, are able to conclude that there must be infinitely many prime numbers.  For those whose Latin is a bit rusty, one of the best modern descriptions of Euler's proof and its consequences appears in chapter VII of \refno{[8]}.

   This short paper uses a simple modification of Euler's argument to obtain new results about the distribution of prime factors of sets of integers.   As this last sentence is something of a mouthful, let us be more clear and specific.  Consider an infinite set $S$ of positive integers such that $\card \lbrace s \in S : s \leq n \rbrace \geq K n^{1/\alpha}$, where $K$ is a constant $> 0$.  Equivalently, let $S$ be the range of a strictly monotone increasing positive integer valued sequence $(s_n)_{n=1}^\infty$, such that $s_n = O(n^\alpha)$.  We shall say that $S$ has polynomial density $\alpha$, even though $\alpha$ may not be an integer.  To avoid trivially degenerate cases, we suppose that $\alpha \geq 1$ and $0 \notin S$.

   Let $P(S)$ be the set of prime factors of elements of $S$, that is let $$P(S) = \lbrace p : p \hbox{ is prime and } p \, | s \hbox{ for at least one } s \in S \rbrace.$$
Let $\pi_S(n) = \card \lbrace p \in P(S) : p \leq n \rbrace$, that is let $\pi_S(n)$ be the number of prime factors $\leq n$ of the elements of $S$.  Clearly $\pi_S(n)$ is the familiar ``prime number function" $\pi(n)$ relativized to $S$.  Our new results are

\theorem{Theorem 1.}

(a) $\sum_{p \in P(S)}  p^{-1/\alpha} = \infty$.

(b) $\sum_{n=1}^\infty  \pi_S(n)/n^{1 + 1/\alpha} = \infty$.

(c) If $\sum_{n=1}^\infty  a_n < \infty$, where $a_n > 0$, then $\pi_S(n)/n^{1 + 1/\alpha} \geq a_n$ infinitely often.

\theorem{Corollary 1.}

(a) There are infinitely many prime factors of the elements of $S$, that is the set $P(S)$ is infinite.

(b) Let $r > 1$.  Then

{\doubleparindent
    $\pi_S(n) \geq n^{1/\alpha}/(\log n)^r$ infinitely often.

    $\pi_S(n) \geq n^{1/\alpha}/\log n (\log\log n)^r$ infinitely often.

    {\sl etc}.\par
} 

\heading{Proof of Corollary 1.}

   Part (a) follows immediately from part (a) of theorem 1.  As for part (b), it is equivalent to

{\doubleparindent
   $\pi_S(n)/n^{1 + 1/\alpha} \geq 1/n\,(\log n)^r$ infinitely often,

   $\pi_S(n)/n^{1 + 1/\alpha} \geq 1/n\,\log n\,(\log\log n)^r$ infinitely often,

   {\sl etc}{.},\par
} 

{\noparindent
and the series $\sum_{n=1}^\infty  1/n\,(\log n)^r$, $\sum_{n=1}^\infty  1/n\,\log n\,(\log\log n)^r$, {\sl etc}{.} are all convergent.  Use theorem 1c.\par}

\heading{Proof of Theorem 1{\rm a}.}

   The proof is a simple modification of Euler's famous proof by infinite series of Euclid's theorem on the infinitude of primes.  Consider the (possibly) infinite product $\prod_{p \in P(S)}  1/(1 - p^{-1/\alpha})$.  We shall show the product diverges.

   To this end, let $N$ be a large positive integer.  Consider the partial product,
$$ \eqalign{\prod_{p \in P(S), p \leq N}  1/(1 - p^{-1/\alpha})
      &= \prod_{p \in P(S), p \leq N}  (1 + p^{-1/\alpha} + p^{-2/\alpha} + \cdots)\cr
      &= \sum_{s \in S, s \leq N}  s^{-1/\alpha} + \hbox{ other positive terms}.}$$
Recall that another way of writing $S$ is $S = \lbrace s_1, s_2, s_3, \ldots \rbrace$, where $(s_j)_{j=1}^\infty$ is a strictly monotone increasing sequence of positive integers such that $s_j \leq K j^{\alpha}$, where $K > 0$ is a constant.  Therefore,
$$\sum_{s \in S, s \leq N}  s^{-1/\alpha} = \sum_{j=1}^n  s_j^{-1/\alpha} \geq \sum_{j=1}^n  K^{-1/\alpha} (j^{\alpha})^{-1/\alpha} = K^{-1/\alpha} \sum_{j=1}^n  j^{-1}.$$
(Here, $n$ is the largest integer such that $s_n \leq N$.)  Since the infinite series $\sum_{j=1}^\infty  j^{-1}$ diverges, so does the infinite product $\prod_{p \in P(S)}  1/(1 -p^{-1/\alpha})$.

   Recall the basic fact that an infinite product $\prod_{n=1}^\infty  (1 + a_n)$, with $a_n \geq 0$, diverges if and only if the infinite series $\sum_{n=1}^\infty  a_n$ diverges.  Thus the divergence of the infinite product
$$\prod_{p \in P(S)}  1/(1 - p^{-1/\alpha})$$
implies the divergence of the infinite series $\sum_{p \in P(S)}\left(1/(1 - p^{-1/\alpha}) - 1\right)$.  The series
$$\sum_{p \in P(S)}\left(1/(1 - p^{-1/\alpha}) - 1\right)$$
is, as we shall show momentarily, term-by-term comparable with the series $\sum_{p \in P(S)}  p^{-1/\alpha}$.  Thus $\sum_{p \in P(S)}  p^{-1/\alpha}$ diverges.  This proves theorem 1a, provided we prove the comparability assertion.

\smallskip

\heading{Proof of the comparability assertion.}

   Note that $p^{-1/\alpha} \leq 2^{-1/\alpha}$ and apply the familiar estimates,
$$ \eqalign{1 + p^{-1/\alpha} &< 1 + p^{-1/\alpha} + p^{-2/\alpha} + \cdots\cr
       &= 1/(1 - p^{-1/\alpha}) = 1 + p^{-1/\alpha}/(1 - p^{-1/\alpha}) \leq 1 + p^{-1/\alpha}/(1 - 2^{-1/\alpha}).}$$
Thus $p^{-1/\alpha} < 1/(1 - p^{-1/\alpha}) - 1 \leq p^{-1/\alpha}/(1 - 2^{-1/\alpha})$, which proves the comparability assertion, and thus proves theorem 1a.

\smallskip

\heading{Proof of Theorem 1{\rm b}.}
$$ \sum_{p \in P(S)}  p^{-1/\alpha} = \lim_{N \rightarrow \infty} \int_1^N t^{-1/\alpha} \, d\pi_S(t),$$
where the right hand side is interpreted as a Stieltjes integral.

Integrate by parts to obtain
$$ \sum_{p \in P(S)}  p^{-1/\alpha} = \lim_{N \rightarrow \infty} \left[ N^{-1/\alpha} \pi_S(N) + (1 + 1/\alpha) \int_1^N \pi_S(t) t^{-(1 + 1/\alpha)} \, dt \right].$$
(The integration by parts is justified because the integrand $t^{-1/\alpha}$ is continuous.)  By theorem 1a, the left hand series diverges to infinity.  We could immediately conclude that
$$ \lim_{N \rightarrow \infty} \int_1^N \pi_S(t) t^{-(1 + 1/\alpha)} \, dt = \infty, $$
if it were not for the troublesome first term on the right, $N^{-1/\alpha} \pi_S(N)$.  To handle this term, divide the argument into two cases.

   In case 1, suppose there is a constant $K > 0$ such that $N^{-1/\alpha} \pi_S(N) < K$ for infinitely many positive integers $N$, say $N_1, N_2, N_3, \ldots$.  Then
$$ \lim_{j \rightarrow \infty} \int_1^{N_j} \pi_S(t) t^{-(1 + 1/\alpha)} \, dt = \infty.$$
Because $\int_1^N \pi_S(t) t^{-(1 + 1/\alpha)} \, dt$ is an increasing function of $N$,
$$ \lim_{N \rightarrow \infty} \int_1^N \pi_S(t) t^{-(1 + 1/\alpha)} \, dt = \infty.$$

   In case 2, suppose there is no such constant $K$.  Then for every $K > 0$, $N^{-1/\alpha} \pi_S(N)$ is eventually $> K$.  Choose any such $K > 0$, then choose $N_K$ so that $N^{-1/\alpha} \pi_S(N) > K$ for all (real) $N > N_K$.  Use the simple estimate,
$$ \int_{N_K}^N \pi_S(t) t^{-(1 + 1/\alpha)} \, dt > K \int_{N_K}^N t^{-1} \, dt$$
to conclude that $\lim_{N \rightarrow \infty} \int_1^N \pi_S(t) t^{-(1 + 1/\alpha)} \, dt = \infty$.

By the integral test (or at least by its proof),
$$ \sum_{n=1}^\infty  \pi_S(n) n^{-(1 + 1/\alpha)} = \infty.$$
This completes the proof of theorem 1b.

\smallskip

\heading{Proof of Theorem 1{\rm c}.}

   Theorem 1c follows immediately from theorem 1b.

\smallskip

\heading{Examples.}

{\noparindent
Example 1.  Polynomial density is needed.  Consider the non-polynomially dense set $S = \lbrace 2^n : n = 1, 2, 3, \ldots \rbrace$.  Then $P(S) = \lbrace 2 \rbrace$.\par}

{\noparindent
Example 2.  The exponent $1/\alpha$ is the best possible.  Let $\alpha$ be an integer $\geq 3$.  We shall construct a set $S$ consisting entirely of primes with polynomial density $\alpha$.  This way, we shall have a good handle on $\pi_S(n)$ and shall easily be able to show $\pi_S(n)$ is asymptotically equal to $n^{1/\alpha}$.\par}

   First, we need to show there is a prime between $n^\alpha$ and $(n + 1)^\alpha$ for all sufficiently large integers $n$.  By a theorem of Iwaniec and Piutz, there is always a prime between $n - n^{23/42}$ and $n$ \cite{5}{p{.}~415}, \refno{[7]}.  If we apply Iwaniec and Piutz's theorem to $n^\alpha$, we see that there is always a prime between $n^{(23/42)\alpha}$ and $n^\alpha$.  But, by the binomial theorem, $(n - 1)^\alpha = n^\alpha - R_\alpha(n)$, where $R_\alpha(n)$ is a polynomial of degree $\alpha - 1$.  Since $(23/42)\alpha < \alpha - 1$ for $\alpha \geq 3$, $n^{(23/42)} < R_\alpha(n)$ for all sufficiently large $n$.  Thus there is always a prime between $(n - 1)^\alpha$ and $n^\alpha$ for all sufficiently large $n$.  Consequently, there is always a prime $p_n$ with $n^\alpha < p_n \leq (n + 1)^\alpha$ for all sufficiently large $n$, say for all $n \geq N_\alpha$.

   Now consider the set $S = \lbrace p_n : n = N_\alpha, N_\alpha + 1, N_\alpha + 2, \ldots \rbrace$.  $S$ has polynomial density $\alpha$, $p_n$ is asymptotically equal to $n^\alpha$, and $\pi_S(n)$ is asymptotically equal to $n^{1/\alpha}$.  This proves the exponent $1/\alpha$ is the best possible in the assertions of theorem 1 and corollary 1b, at least when $\alpha$ is an integer $\geq 3$.

   As we shall see shortly, this remains true when $\alpha = 1$.  Furthermore, the exponent $1/\alpha$ is {\sl probably\/} the best possible for $\alpha = 2$, as it is widely believed, though still unproven, that for every integer $n \geq 2$, there is a prime $p_n$ between $n^2$ and $(n + 1)^2$ \cite{5}{p{.}~19}.

{\noparindent
Example 3.  The $\log n$ denominator (though possibly not the exponent $r$) is needed.  Let $S$ be the set of positive integers.  Then $S$ has polynomial density 1, and $\pi_S(n) = \pi(n)$.  By Tschebyshev's inequalities,
$$ m n/\log n \leq \pi(n) \leq M n/\log n $$
for all $n \geq 2$, where $m$ and $M$ are suitably chosen constants.\par}

  This example also shows that the exponent $1/\alpha$ is the best possible for $\alpha = 1$.

{\noparindent
Example 4.  Furstenberg's proof and Gotchev's theorem.  Again let $S$ be the set of positive integers.  Furstenberg (and later Golomb) gave elegant {\sl topological proofs\/} of Euclid's theorem that $P(S)$ is infinite \refno{[2]}, \refno{[3]}.  This work has been extended in various directions, for example to the setting of abstract ideal theory (see  \refno{[9]} and \refno{[13]}).  Furstenberg's proof is also provides an important beginning example in the theory of profinite groups \cite{11}{p{.}~476}\par}

   Gotchev recently extended Furstenberg's work and gave a {\sl topological proof\/} that $P(S)$ is infinite if $S$ is the integer range of a polynomial $Q$ of degree $\alpha \geq 1$ with integer coefficients \refno{[4]}.  In other words, if $S = \lbrace |Q(n)| : n = 1, 2, 3, \ldots \rbrace$, Gotchev proved by {\sl topological means\/} that $P(S)$ is infinite.  As $S$ is a set of polynomial density $\alpha$, corollary 1a is Gotchev's theorem in this case, and corollary 1b is a quantitative version of Gotchev's theorem.

   We should note that Gotchev's theorem seems to be a folk theorem, probably rediscovered many times before Gotchev stated it.  But the {\sl topological proof\/} is definitely due to Gotchev.

{\noparindent
Example 5.  Sets of polynomial density $\alpha$ are much more general than ranges of polynomials.  In fact, they may even be non-computable.  For example, let $T(n), n = 1, 2, 3, \ldots$ be your favorite (computable) enumeration of Turing machines, let $k(n) = 1$ if $T(n)$ does not halt, and let $k(n) = 2$ otherwise.  Let $S = \lbrace 3n + k(n) : n = 1, 2, 3, \ldots \rbrace$.  Then $S$ is a non-computable set of polynomial density 1.\par}

   Although we know practically nothing about $S$ because it is not computable, we do know by corollary 1a that its set of prime factors is infinite, and we do have quantitative information about the distribution $P(S)$ by theorem 1.

{\noparindent
Example 6.  One might think that the complement of a set too sparse to be of polynomial density $\alpha$ for any $\alpha$ must be a set of polynomial density 1.  An interesting counter example due to Joshua Zelinsky \refno{[15]} shows that this is not true.  Zelinsky's idea is to consider sets with very large gaps, that is highly lacunary sets.\par}

   To be specific, let $S_n = \lbrace n : 2^{2^n} \leq n < 2^{2^{n + 1}} \rbrace$, and let $S = \bigcup_{n \hbox{\smrm \ odd}} S_n$.  Then $S$ has a gap between $2^{2^n}$ and $2^{2^{n + 1}} - 1$ for each even $n$,  and the complement of $S$ has a gap between $2^{2^n}$ and $2^{2^{n + 1}} - 1$ for each odd $n$   Thus neither $S$ nor its complement is of polynomial density $\alpha$ for any $\alpha$.

\smallskip

\heading{Remarks.}

{\noparindent
Remark 1. Dirichlet's theorem on primes in arithmetic progressions states that if $a$ and $b$ are relatively prime integers, and $S = \lbrace an + b : n = 1, 2, 3, \ldots \rbrace$, then $S$ {\sl contains\/} infinitely many primes \cite{8}{chapter VII}.  This is much stronger than the assertion of Gotchev's theorem, and corollary 1a, that $P(S)$ is infinite.  In general, it seems to be much more difficult to show that a set {\sl contains\/} infinitely many primes than it is to show that a set has infinitely many prime factors.\par}

   An extremely interesting challenge would be to give a simple, fairly elementary proof of Dirichlet's theorem, perhaps using the ideas of this paper, or perhaps using the topological ideas of Furstenberg, Golomb and Gotchev \cite{11}{p{.}~476}.  It turns out that if one could show that an arithmetic progression with $a$ and $b$ relatively prime contains {\sl one\/} prime, then it contains infinitely many primes \refno{[6]}, \cite{14}{p{.}~124}.

{\noparindent
Remark 2.  Buniakowski conjectured in 1857 that if $Q$ is a polynomial with integer coefficients, if $S = \lbrace |Q(n)| : n = 1, 2, 3, \ldots \rbrace$, and if there is no common factor of all the elements of $S$, then $S$ contains infinitely many primes \refno{[1]}.  The Buniakowski conjecture is one of the great unsolved questions of number theory \cite{10}{p{.}~323}, \cite{12}{p{.}~452}, see also \cite{14}{p{.}~127~ff{.}}.  (In fact, it is even unknown if there are infinitely many primes of the form $n^2 + 1$ \cite{5}{p{.}~19}.) Once again, this illustrates how assertions about the {\sl existence of primes\/} in sets seem to be far more deeper that assertions about the distribution of {\sl prime factors}.\par}

\smallskip

\heading{Open Questions.}

{\noparindent
Question 1.  Is there an example for $\alpha > 1$ where both the exponent $1/\alpha$ and the denominator $\log n$ appear at the same time?  Of course, there is such an example for $\alpha = 1$; simply let $S$ be the set of positive integers and use Tschebyshev's inequalities.\par}

{\noparindent
Question 2.  There cannot be a prime number theorem for $P(S)$, or even a pair of Tschebyshev inequalities (see example 2).  But can corollary 1b be improved to a {\sl one-sided\/} Tschebyshev inequality?  In other words, for every set $S$ of polynomial density $\alpha$, is it true that
$$ m_S \, n^{1/\alpha}/\log n \leq \pi_S(n),$$
for all $n \geq 2$ and some suitably chosen constant $m_S > 0$ (depending on $S$)?\par}

\endsmallparskip


\smallskip
\heading{References.}
\smallskip

{\noparindent

\refno{1{.}}~V.~Buniakowski, Sur les diviseurs num\'eriques invariables des fonctions rationnelles entri\`eres, {\sl Acad. Sci. P\'etersbourg Mem.} {\bf 6} {\sl sci. math. et phys.} {\bf 6} (1857), 305-329.

{\refno{2{.}}~H{.}~Furstenberg, On the infinitude of primes, {\sl Amer{.} Math{.} Monthly\/} {\bf 62} (1955), 353.

\refno{3{.}}~S{.}~W{.}~Golomb, A connected topology for the integers, {\sl Amer{.} Math{.} Monthly\/} {\bf 66} (1959), 663-665.

\refno{4{.}}~I{.}~Gotchev, On the Infinitude of Prime Divisors of the Range of a Non-constant Polynomial, to appear.

\refno{5{.}}~G{.}~H{.}~Hardy and E{.}~M{.}~Wright, {\sl An Introduction to the Theory of Numbers}, 5th ed{.}, Oxford University Press, Oxford, 1993.

\refno{6{.}}~V{.}~Hanly, Solution to problem E1218, {\sl Amer{.} Math{.} Monthly\/} {\bf 64} (1957), 742.

\refno{7{.}}~Iwaniec and Piutz, {\sl Monatshefte f{.} Math}{.} {\bf 98} (1984), 115-143.

\refno{8{.}}~A.~Knapp, {\sl Elliptic Curves}, Princeton University Press, Princeton, NJ, 1992.

\refno{9{.}}~J{.}~Knopfmacher and S{.}~Porubsky, Topologies related to arithmetical properties of integral domains, {\sl Expositiones Mathematicae\/} {\bf 15}(2) (1997), 131-148.

\refno{10{.}}~S{.}~Lang, {\sl Algebra}, 3rd ed{.}, Addison Wesley, Reading, MA, 1993.

\refno{11{.}}~A{.}~Lubotzky (Reviewer), {\sl Profinite groups\/} by J{.}~S{.} Wilson, {\sl Profinite groups\/} by L{.} Ribes and P{.} Zalesskii, {\sl Analytic pro-p groups\/} by J{.}~D{.}~Dixon, M{.}~P{.}~F{.} du Satoy, A{.}~Mann and D{.}~Segal, {\sl New horizons in pro-p groups\/} by M{.}~du Satoy, D{.}~Segal and A{.} Shalev (eds.), {\sl Bull{.} Amer{.} Math{.} Soc}{.} (N.S.) {\bf 38}(4) (2001), 475-479.  Erratum, {\bf 39}(2) (2002), 299.

\refno{12{.}}~M{.}~Ram Murty, Prime numbers and irreducible polynomials, {\sl Amer{.} Math{.} Monthly\/} {\bf 109}(5) (2002), 452-458.

\refno{13{.}}~S{.}~Porubsky, Arithmetically related ideal topologies and the infinitude of primes, {\sl Questiones Mathematicae\/} {\bf 24} (2001).

\refno{14{.}}~W{.}~Sierpinsky, {\sl Elementary Theory of Numbers}, Warszawa, Panstwowe Wydawn{.} Naukowe, 1964.

\refno{15{.}}~J{.}~Zelinsky, private communication.

} 

\bye